# Regression with Strongly Correlated Data

Christopher S. Jones[(1)], John M. Finn[(1)] and Nicolas Hengartner[(2)]
[(1)]T-15, Plasma Theory and [(2)]D-1, Statistical Sciences
Los Alamos National Laboratory, Los Alamos, NM 87545

1st September 2006

**Abstract**

*This paper discusses linear regression of strongly correlated data that arises, for example, in magnetohydrodynamic equilibrium reconstructions. We have proved that, generically, the covariance matrix of the estimated regression parameters for fixed sample size goes to zero as the correlations become unity. That is, in this limit the estimated parameters are known with perfect accuracy. Simple examples are shown to illustrate this effect and the nature of the exceptional cases in which the estimate covariance does not go to zero.*

## 1  Introduction

Magnetohydrodynamic (MHD) equilibrium reconstructions play a vital role in the analysis of the states of plasmas in magnetic confinement devices such as tokamaks [1, 2, 3, 4]. Typically, such reconstructions are performed by least squares fitting of the nonlinear Grad-Shafranov equation to measurements of the magnetic field at spatially distinct points on the boundary of the device, complemented by measurements of the interior conditions of the plasma. Dynamical fluctuations associated with plasma turbulence occurring on relatively short time scales are modeled as stochastic noise in the reconstructions, and these fluctuations may exhibit strong correlations in space and time.

Operationally, regression with correlated errors is understood [5]. However, the least squares equilibrium reconstruction studies of reference [6] exhibited unexpected properties. As a function of the degree of correlation (discussed further below), the variance of the fitted parameters was observed to have a maximum. Past this peak, the estimate covariance matrix converged to zero as the simulated measurements became fully correlated. The objective of this paper is to explore this phenomenon and show how generic it is. While our analysis is restricted to linear regression because the effect is most transparent there, our analysis and results are readily extended to nonlinear regression and were indeed first observed in a nonlinear context [6]. Specifically, consider the linear model

$$Y_i = X_i^t \beta + \eta_i, \quad i = 1, 2, \ldots, n \tag{1}$$

where the covariates and the parameters $X_i, \beta \in \mathbb{R}^m$ with $m \leq n$, and the disturbances $\eta_i$ have mean zero, variances $\mathbb{E}[\eta_i^2] = \sigma_i^2$ and covariances $\mathbb{E}[\eta_i \eta_j] = \Sigma_{ij} = \sigma_i \sigma_j \varrho_{ij}$. For our analysis and discussion, it is convenient to rewrite (1) in vector form

$$\mathsf{Y} = \mathsf{X}\beta + \eta, \tag{2}$$

where the vectors $\mathsf{Y}, \eta \in \mathbb{R}^n$, $\mathsf{X}$ is an $n \times m$ design matrix, $\mathbb{E}[\eta] = 0$ and $\mathbb{E}[\eta \eta^t] = \Sigma = \mathsf{SRS}$, where the deviations are $\mathsf{S} = \mathrm{diag}(\sigma_1, \sigma_2, \ldots, \sigma_n)$ and $\mathsf{R}$ is the matrix of correlation coefficients $\varrho_{ij}$. To simplify our exposition, we shall begin by assuming that the design matrix $\mathsf{X}$ is of rank $m$, and unless otherwise stated, further assume that the covariance matrix $\Sigma$ is of full rank $n$.

It is well known that the best linear unbiased estimators (BLUE) $\hat{\beta}$ for $\beta$ is the minimizer (with respect to $\beta$) of the quadratic

$$\chi^2(\beta) = (\mathsf{Y} - \mathsf{X}\beta)^t \Sigma^{-1} (\mathsf{Y} - \mathsf{X}\beta). \tag{3}$$

Under the stated assumptions, the BLUE $\hat{\beta}$ is then given by the unique solution of the normal equations

$$(\mathsf{X}^t \Sigma^{-1} \mathsf{X}) \hat{\beta} = \mathsf{X}^t \Sigma^{-1} \mathsf{Y}, \tag{4}$$

and the covariance matrix of the estimators $\hat{\beta}$ is given by

$$\mathsf{V} = (\mathsf{X}^t \Sigma^{-1} \mathsf{X})^{-1}. \tag{5}$$

Our main result is that, if the vector of signed standard deviations $(e_1 \sigma_1, \ldots, e_n \sigma_n)^t$ of the disturbances, with signs $e_j = \pm 1$ determined consistently with respect to the correlations $\varrho_{ij}$ (as explained in Section 3), does not lie in the column space of the design matrix X, then the sum of the variances, trace(V), of the estimated regression coefficients converges to zero as the noise becomes fully correlated. Full correlation is defined as the limit

$$\min_{i,j} \frac{|\Sigma_{ij}|}{\sigma_i \sigma_j} = \min_{i,j} |\varrho_{ij}| \longrightarrow 1.$$

This paper is organized as follows. Section 2 sets the stage by introducing the simplest possible model exhibiting our main result. This simple model is related to Peelle's Pertinent Puzzle [7, 8, 9] and motivates the main result, which is presented in Section 3. Section 4 deals with the issue of increasing the number of measurements in physical problems such as MHD equilibrium reconstruction, by adding measurements at necessarily more closely packed spatial positions. In Section 5 we deal with various generalizations involving the rank of the design matrix and the data covariance matrix. A summary and conclusions are presented in Section 6.

## 2 A simple example

This section presents a simple example consisting of estimating the common mean from a pair of highly correlated observations. Specifically, we wish to estimate $\mu$ from a pair of observations

$$\begin{pmatrix} y_1 \\ y_2 \end{pmatrix} = \begin{pmatrix} 1 \\ 1 \end{pmatrix} \mu + \begin{pmatrix} \eta_1 \\ \eta_2 \end{pmatrix} = \mathsf{X}\mu + \begin{pmatrix} \eta_1 \\ \eta_2 \end{pmatrix}, \qquad (6)$$

where the disturbances $(\eta_1, \eta_2)$ have mean zero and covariance and precision matrices given by

$$\Sigma = \begin{pmatrix} \sigma_1^2 & \varrho\sigma_1\sigma_2 \\ \varrho\sigma_1\sigma_2 & \sigma_2^2 \end{pmatrix} \quad \text{and} \quad \Sigma^{-1} = \frac{1}{1-\varrho^2} \begin{pmatrix} \tau_1^2 & -\varrho\tau_1\tau_2 \\ -\varrho\tau_1\tau_2 & \tau_2^2 \end{pmatrix},$$

with $\tau_1 = 1/\sigma_1$, $\tau_2 = 1/\sigma_2$. Direct calculation shows that the weighted least squares estimate that satisfies the normal equations (4) is

$$\begin{aligned} \hat{\mu} &= \frac{(\tau_1^2 - \varrho\tau_1\tau_2) y_1 + (\tau_2^2 - \varrho\tau_1\tau_2) y_2}{\tau_1^2 - 2\varrho\tau_1\tau_2 + \tau_2^2} \\ &= \frac{(\tau_1^2 - \varrho\tau_1\tau_2) y_1 + (\tau_2^2 - \varrho\tau_1\tau_2) y_2}{(1-\varrho)(\tau_1^2 + \tau_2^2) + \varrho(\tau_1 - \tau_2)^2} \\ &= w_1(\sigma_1, \sigma_2, \varrho) y_1 + w_2(\sigma_1, \sigma_2, \varrho) y_2, \end{aligned} \qquad (7)$$

where we have written the estimate in terms of weights $w_1$ and $w_2$ on the last line. The variance of this estimate, computed from equation (5), is

$$\begin{aligned} V(\hat{\mu}) &= \frac{1-\varrho^2}{\tau_1^2 - 2\varrho\tau_1\tau_2 + \tau_2^2} \\ &= \frac{1-\varrho^2}{(1-\varrho)(\tau_1^2 + \tau_2^2) + \varrho(\tau_1 - \tau_2)^2}. \end{aligned} \quad (8)$$

The behavior of equation (8) as a function of $\varrho$ for $\sigma_1 = 1$ and various $\sigma_2$ is shown in Fig. 1, showing clearly the limits $\varrho \to \pm 1$. In this simple example (and as displayed in Fig. 1), the estimate variance goes to zero as $\varrho \to \pm 1$, with the exception of $\varrho \to 1$ with $\sigma_1 = \sigma_2$. Notice that, for the exceptional case, the vector $(\sigma_1, \varrho\sigma_2) \to (\sigma_1, \sigma_1)$ is in the range space of the design matrix $\mathsf{X}$ in equation (6).

Regarding the exceptional case of $\sigma_1 = \sigma_2$, the estimate for $\mu$ for all $\varrho \in [-1, 1]$ is the sample mean

$$\hat{\mu} = \frac{y_1 + y_2}{2}, \quad \text{which has variance} \quad V(\hat{\mu}) = \frac{1+\varrho}{2}\sigma^2. \quad (9)$$

The variance vanishes for $\varrho \to -1$ while for $\varrho \to 1$, it is

$$\lim_{\varrho \to 1} V(\hat{\mu}) = \sigma^2. \quad (10)$$

The interpretation of equation (10) is the following. For a single measurement $y_1 = \mu + \eta_1$, the variance of $\hat{\mu} = y_1$ equals the variance of $y_1$, namely $\sigma^2$. For $\varrho \to 1$ and $\sigma_1 = \sigma_2$, the second measurement $y_2$ is equal to $y_1$ and provides no new information, leaving $V = \sigma^2$. This is the familiar interpretation of positive correlations, which are often assumed to imply redundancy in measurement. We shall see below, however, that in the generic case positive correlations may provide leverage with which to determine more accurately the estimator.

In the generic case of $\sigma_1 \neq \sigma_2$, the possibility arises of *negative weighting*, in which one of the two weights, $w_1$ or $w_2$, becomes negative. In particular, from equation (7) and assuming without loss of generality $\sigma_1 > \sigma_2$, $y_1$ is weighted negatively ($w_1 < 0$) if $\varrho > \sigma_2/\sigma_1$. Also, from equation (8) we conclude that $V(\hat{\mu})$ has a maximum with respect to $\varrho$ at $\varrho = \sigma_2/\sigma_1$. At the maximum variance, note that the estimate (7) becomes $\hat{\mu} = y_2$ and equation (8) gives $V(\hat{\mu}) = \sigma_2^2$, i.e. $\hat{\mu}$ becomes independent of the measurement $y_1$. For $\varrho > \sigma_2/\sigma_1$, $V(\hat{\mu})$ decreases with respect to $\varrho$. The possibility of negative weighting, leading to an estimate outside the range of the measurements, has been noted with surprise by the nuclear data community, where it is known as Peelle's pertinent puzzle[7, 8, 9]. These investigations did not, however, remark on the observation that for fixed $\sigma_1 \neq \sigma_2$, the appearance of negative weighting (as $\varrho$ is increased) coincides with the decrease of the variance $V(\hat{\mu})$.

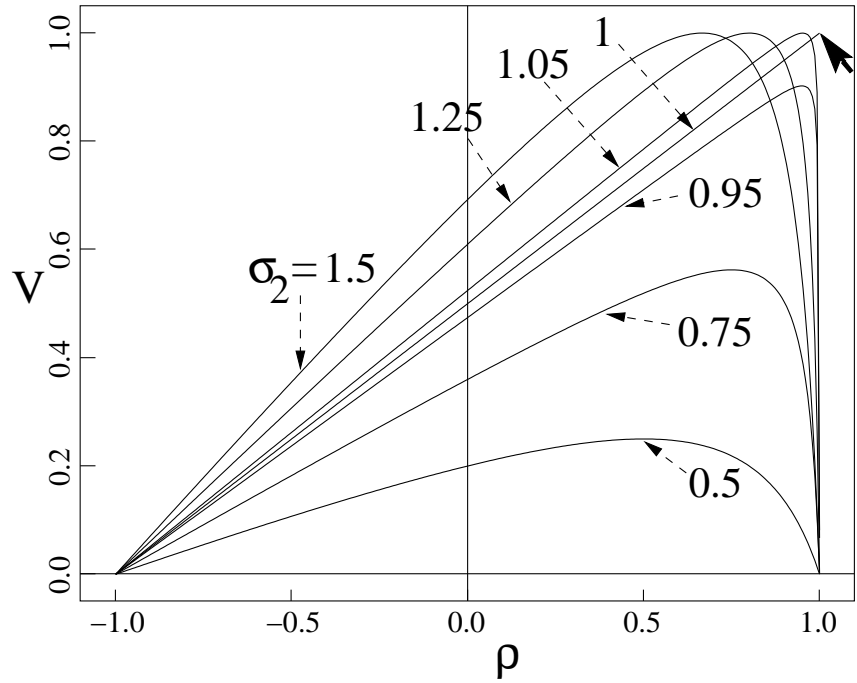

Figure 1: Estimate variance $V$ as a function of the correlation coefficient $\varrho$ for $\sigma_1 = 1$ and seven values of $\sigma_2$ (0.5, 0.75, 0.95, 1, 1.05, 1.25, 1.5). For all cases, $V(\hat{\mu}) \to 0$ as $\varrho \to -1$. For $\sigma_1 \neq \sigma_2$, $V(\hat{\mu}) \to 0$ as $\varrho \to 1$. For $\sigma_2 = \sigma_1 = 1$, $V(\hat{\mu}) = \sigma_1^2 = 1$ at $\varrho = 1$ (large arrow).

In the limit of full correlation ($\varrho \to 1$), the estimator for $\mu$ in the generic case of $\sigma_1 \neq \sigma_2$ is

$$\hat{\mu} = (\tau_1 y_1 - \tau_2 y_2) / (\tau_1 - \tau_2), \tag{11}$$

with negative weighting on $y_1$ if $\tau_1 < \tau_2$ ($\sigma_1 > \sigma_2$). This estimate has the variance (in the same limit)

$$V(\hat{\mu}) \to \frac{2(1-\varrho)}{(\tau_1 - \tau_2)^2} \to 0, \tag{12}$$

To interpret these results, consider Fig. 2a. The noise contributions $\eta_1$ and $\eta_2$ will necessarily have the same sign if $\varrho \to 1$. If, for example, they are both positive, then $y_1$ and $y_2$ will both be above $\mu$, and an unbiased estimate will be possible only with negative weighting as in equation (11). To be more specific, for $\varrho \to 1$ but $\sigma_1 > \sigma_2$, we will have $\eta_2 = \eta_1 \sigma_2 / \sigma_1$ or

$$y_1 = \mu + \sigma_1 \alpha, \quad y_2 = \mu + \sigma_2 \alpha, \tag{13}$$

where $\alpha$ is a single random variable with zero mean and unit variance. The BLUE (11) chooses the correct (negative) weighting to give the exact result $\hat{\mu}$ using a *single* realization $\alpha$ of the noise. The noise term can be eliminated from equation (13), giving

$$\frac{y_1}{\sigma_1} - \frac{y_2}{\sigma_2} = \left(\frac{1}{\sigma_1} - \frac{1}{\sigma_2}\right)\mu, \tag{14}$$

in agreement with the estimate in equation (11). The variance of this estimate is zero because the noise has been eliminated. Equivalently, a different realization of the noise (different $\alpha$) yields the same result. If, on the other hand, $\sigma_1 = \sigma_2$, then the measurements $y_1$ and $y_2$ are identical and the process leading to equation (14) cannot be followed.

An alternate geometric interpretation is shown in Fig. 2b. In this case, for normally distributed noise, the measurements are distributed according to a probability density function proportional to $\exp(-\chi^2/2)$. The level sets of this function (contours of $\chi^2$) are ellipses which circumscribe regions within which the measurements may be found with a given probability. For $\varrho$ close to unity, the measurements are expected to be found within a thin ellipse whose major axis has slope $\sigma_2/\sigma_1$. For $\varrho \to 1$, the two eigenvalues of the covariance matrix $\Sigma$ are $\lambda_1 = \sigma_1^2 + \sigma_2^2$ (trace) and $\lambda_2 = 0$, and the ellipses become infinitely thin, i.e. line segments. In this limit the estimate must be on the intersection of the line $y_2 - \mu = \sigma_2/\sigma_1 (y_1 - \mu)$ and the axis $y_1 = y_2$, giving $\hat{\mu} = \mu$ with zero uncertainty.

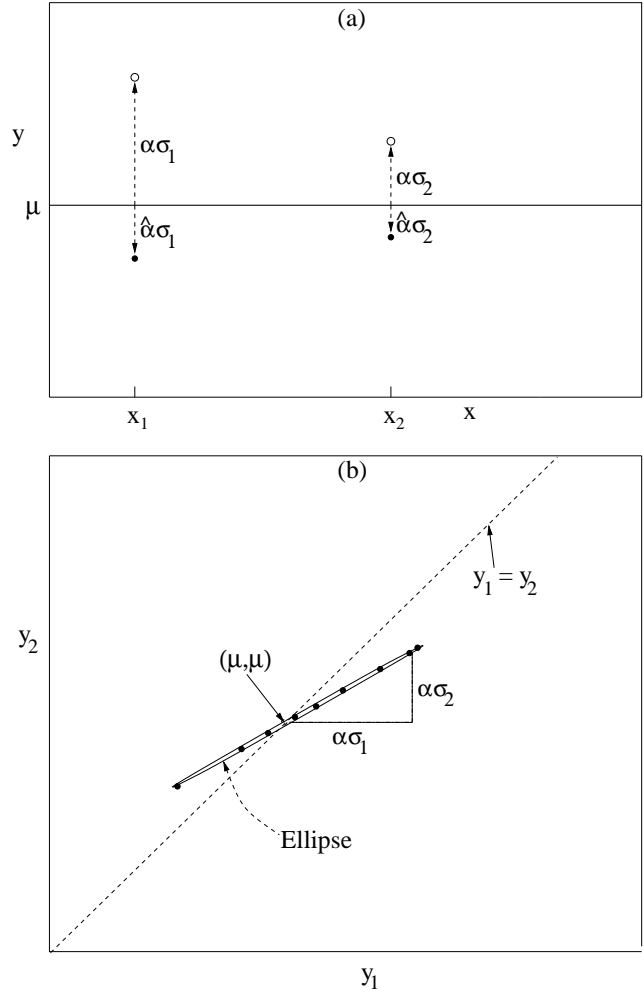

Figure 2: Estimating a constant in the limit of large correlations $\varrho \to 1$. Two noisy measurements of $\mu$ are made (a) at positions $x = x_1$, $x = x_2$ (open circles), corresponding to one value of $\alpha$. The estimate, which has zero variance, has $(y_1 - \mu)/\sigma_1 = (y_2 - \mu)/\sigma_2$, exactly determining $\hat{\mu}$ if $\sigma_1 \neq \sigma_2$. A second set of correlated measurements (solid circles), with a second value $\hat{\alpha}$, leads to the same estimate. In (b), for $\sigma_1 > \sigma_2$ and $\varrho \to 1$, when the ellipse at constant $\chi^2$ collapses, each measurement point $(y_1, y_2)$ has $y_1 = \mu + \alpha\sigma_1$, $y_2 = \mu + \alpha\sigma_2$ for some $\alpha$. Since the estimate must also be along the line $y_1 = y_2$, the estimate for $\varrho \to 1$ gives $\hat{\mu} = \mu$ with zero uncertainty.

As we shall see in the next section, this phenomenon is the rule rather than the exception, and it is linked to the emergence of a *noise-free subspace* in the limit of large correlations.

## 3  Regression in the limit of full correlation

This section generalizes the example of Section 2 and shows under what conditions the variance of the weighted least squares regression estimator for the linear model (2) converges to zero in the limit of full correlation. We consider classes of covariances for the disturbances $\eta$ for which we may define the parameter

$$\kappa = \max_{i,j} \left( 1 - \frac{|\Sigma_{ij}|}{\sigma_i \sigma_j} \right) = \max_{i,j} \left( 1 - |\varrho_{ij}| \right)$$

The covariance matrix can be written as $\Sigma = \mathsf{SRS}$, where $\mathsf{S} = \text{diag}(\sigma_1, \sigma_2, \ldots, \sigma_n)$ and $\mathsf{R}$ is the matrix of correlation coefficients $\varrho_{ij}$. In the limit of $\kappa = 0$, the correlation matrix is the rank one matrix $\mathsf{R} = ee^t$, where $e_j = \pm 1$. A simple example of such a class of covariance matrices is the autocorrelation model. In this model, the disturbances $\eta$ consist of a mean zero random vector with covariance matrix $\Sigma = \mathbb{E}[\eta \eta^t] = \mathsf{SRS}$, where the correlation matrix is

$$\mathsf{R} = \begin{pmatrix} 1 & \varrho & \ldots & \varrho^{n-1} \\ \varrho & 1 & \ldots & \varrho^{n-2} \\ \ldots & \ldots & \ddots & \ldots \\ \varrho^{n-1} & \varrho^{n-2} & \ldots & 1 \end{pmatrix}, \tag{15}$$

or $\varrho_{ij} = \varrho^{|i-j|}$. For this model, we may identify $\kappa \equiv (1 - |\varrho|^{n-1})$. We are interested in the behavior of the variance $(\mathsf{X}^t \Sigma^{-1} \mathsf{X})^{-1}$ of the BLUE $\hat{\beta}_\kappa = (\mathsf{X}^t \Sigma^{-1} \mathsf{X})^{-1} \mathsf{X}^t \Sigma^{-1} \mathsf{Y}$ as $\kappa \to 0$. In this section, we emphasize the dependence of the estimated regression parameter on the parameter $\kappa$ by subscripting the estimate $\hat{\beta}_\kappa$. Note that in the limit $\kappa = 0$, the correlation matrix of the autocorrelation model (15) is $\mathsf{R} = ee^t$, where $e$ has entries $e_j = 1$ for $\varrho \to 1$, or $e_j = (-1)^j$ for $\varrho \to -1$.

We begin by presenting a heuristic that shows how vanishing of the estimator variance is related to the emergence of a noise-free subspace in the limit of large correlations. A more rigorous proof is then presented at the end of this section.

Fix the sample size to $n > m$, and let $\mathsf{Q} = [v_1, v_2, \ldots, v_n]$ denote the matrix of normalized column eigenvectors of the covariance matrix $\Sigma$, corresponding to the eigenvalues $\lambda_1 > \lambda_2 > \ldots > \lambda_n > 0$. As $\kappa \longrightarrow 0$, we have $\Sigma_{ij} \to \sigma_i e_i \sigma_j e_j$, so that

$$\lambda_1 = \sum_{j=1}^{n} \sigma_j^2 + o(1), \ \ \lambda_k = o(1), \ k \geq 2, \ \text{and} \ \ v_{1j} = e_j \frac{\sigma_j}{\sqrt{\lambda_1}} + o(1). \quad (16)$$

Consider the transformation

$$\mathsf{Z} = \mathsf{Q}^t \mathsf{Y} = (\mathsf{Q}^t \mathsf{X})\beta + \mathsf{Q}^t \eta = \tilde{\mathsf{X}}\beta + \Lambda^{1/2}\xi, \quad (17)$$

where $\Lambda = \text{diag}(\lambda_1, \ldots, \lambda_n)$ and $\xi$ is a vector of uncorrelated disturbances with mean zero and unit variance.

In light of (16), $\lambda_j = 0, j \geq 2$ in the limit of $\kappa = 0$, so that the transformed variables $Z_j$ ($j \geq 2$) are noise-free. That is, the subset of equations (17) with $j \geq 2$ are not subjected to noise in the limit $\kappa = 0$:

$$Z_j = v_j^t \mathsf{X}\beta, \ \ j = 2, \ldots, n \quad (18)$$

Further, if $v_1$ does not lie in the column space (the range space) of $\mathsf{X}$, then the system of equations (18) has a unique solution $\beta^\star$ (recall $n > m$ and $\mathsf{X}$ is assumed of full rank). We call the space spanned by the eigenvectors $v_j$ for $j \geq 2$ the *noise-free subspace*. Since the weighted least squares estimator $\hat{\beta}_\kappa$ has the smallest variance among all linear estimators for $\beta$, we conclude that, assuming that the BLUE exists in the limit,

$$\text{trace}(\mathsf{V}(\hat{\beta}_{\kappa=0})) \leq \text{trace}(\mathsf{V}(\beta^\star)) = 0$$

for every fixed sample size $n > m$, suggesting that $\hat{\beta}_{\kappa=0} = \beta^\star$.

Conversely, suppose that $v_1$ lies in the column space of $\mathsf{X}$. Let us reparametrize the column space of $\mathsf{X}$ via a linear one-to-one transformation $\mathsf{W}$ such that the first column of $(\mathsf{X}\mathsf{W})$ is $v_1$ and all others columns are orthogonal to $v_1$. In this way we identify the linear combination of elements of $\beta$ which are subjected to noise. Setting $\gamma = \mathsf{W}^{-1}\beta$, we get

$$\mathsf{Y} = \mathsf{X}\mathsf{W}\mathsf{W}^{-1}\beta + \eta = (\mathsf{X}\mathsf{W})\gamma + \eta.$$

Then (16) implies that in the limit of $\kappa = 0$,

$$Z_1 = v_1^t \mathsf{Y} = \gamma_1 + \left(\sum_{j=1}^{n} \sigma_j^2\right)^{1/2} \xi_i, \tag{19}$$

$$Z_j = v_j^t \mathsf{Y} = v_j^t (\mathsf{XW})\gamma, \quad j = 2, \ldots, n. \tag{20}$$

Because the first column of $(\mathsf{XW})$ is $v_1$, equations (20) are independent of $\gamma_1$ and determine exactly the parameters $\gamma_2, \ldots, \gamma_m$. The estimator of $\gamma_1$, a particular linear combination of the elements of $\beta$, is subjected the total noise of the system (i.e. the variance of the estimate $\hat{\gamma}_1$ is the sum of the variances of the original measurements). Hence in the limit as $\kappa \longrightarrow 0$, the total variance for the BLUE for $\gamma$ is

$$\text{trace}(\mathsf{V}(\hat{\gamma}_\kappa)) = \sum_{j=1}^{n} \lambda_j = \sum_{j=1}^{n} \sigma_j^2 > 0,$$

which implies that $\text{trace}(\mathsf{V}(\hat{\beta}_\kappa)) > 0$.

This heuristic identifies, in the limit of $\kappa = 0$, the subspace orthogonal to $v_1$ as a noise-free subspace that enables perfect estimation of $\beta$. However, this argument does not prove that the variance of the BLUE for $\beta$ converges to zero with $\kappa \longrightarrow 0$, because our argument lets $\kappa$ converge to zero first before estimating $\beta$ and showing that it resulted in an estimate that had zero variance. Theorem 1 below gives a rigorous proof of our claim.

THEOREM 1. *For fixed sample size $n > m$, suppose that the design matrix $\mathsf{X}$ is of rank $m$. If the eigenvector $v_1$ associated with the largest eigenvalue of the limiting covariance matrix (when $\kappa$ goes to zero) does not lie in the column space of the design matrix $\mathsf{X}$, then the total variance of the least squares estimate approaches zero in the limit as $\kappa \longrightarrow 0$.*

**Proof.** Denote by $\mathsf{Q} = [v_1, \ldots, v_n]$ the matrix of column eigenvectors associated with the eigenvalues $\lambda_1 \geq \lambda_2 \geq \ldots \geq \lambda_n > 0$ of the covariance matrix $\Sigma$. As before, the continuity of the eigenvectors and eigenvalues as a function of $\Sigma$, implies that as $\kappa \longrightarrow 0$, the eigenvalues and first eigenvector behave as in (16). We again consider the transformation defined in equation (17). Suppose that $v_1$ is not in the range of the design matrix $\mathsf{X}$. Then

$$Z_1 = v_1^t \mathsf{X}\beta + \sqrt{\lambda_1}\xi_1 = v_1^t \mathsf{X}\beta + \sqrt{\sigma_1^2 + \ldots + \sigma_n^2}\xi_1,$$

so that $Z_1$ is a linear combination of the regression parameters $\beta_i$ containing all the noise. Then the BLUE for $\beta$ is

$$\hat{\beta}_\kappa = (\tilde{\mathsf{X}}^t_{-1}\Lambda^{-1}_{-1}\tilde{\mathsf{X}}_{-1})^{-1}\tilde{\mathsf{X}}^t_{-1}\Lambda^{-1}_{-1}Z_{-1},$$

where $Z_{-1}$ is the vector $(Z_2,\ldots,Z_n)^t$, $\Lambda_{-1} = \mathrm{diag}(\lambda_2,\ldots,\lambda_n)$, and $\tilde{\mathsf{X}}_{-1} = [v_2, v_3, \ldots, v_n]^t\mathsf{X}$. The covariance matrix of $\hat{\beta}$ is

$$\mathsf{V}(\hat{\beta}_\kappa) = (\tilde{\mathsf{X}}^t_{-1}\Lambda^{-1}_{-1}\tilde{\mathsf{X}}_{-1})^{-1},$$

so that the largest eigenvalue satisfies

$$\begin{aligned}
\sup_{\|a\|=1} a^t\mathsf{V}(\hat{\beta}_\kappa)a &= \sup_{\|a\|=1} a^t(\tilde{\mathsf{X}}^t_{-1}\Lambda^{-1}_{-1}\tilde{\mathsf{X}}_{-1})^{-1}a \\
&= \left[\inf_{\|a\|=1} a^t\tilde{\mathsf{X}}^t_{-1}\Lambda^{-1}_{-1}\tilde{\mathsf{X}}_{-1}a\right]^{-1} \\
&= \left[\inf_{\forall a} \frac{\left(\tilde{\mathsf{X}}_{-1}a\right)^t\Lambda^{-1}\left(\tilde{\mathsf{X}}_{-1}a\right)}{\|\tilde{\mathsf{X}}_{-1}a\|^2}\frac{\|\tilde{\mathsf{X}}_{-1}a\|^2}{\|a\|^2}\right]^{-1}.
\end{aligned}$$

Noting that, for strictly positive quantities, the infimum of a product is greater than or equal to the product of infima, we have

$$\sup_{\|a\|=1} a^t\mathsf{V}(\hat{\beta}_\kappa)a \leq \lambda_2 / \inf_{\|a\|=1} \|\tilde{\mathsf{X}}_{-1}a\|^2.$$

In light of (16) and the full rank of $\tilde{\mathsf{X}}_{-1}$, the latter converges to zero with $\kappa \longrightarrow 0$. □

**Remark:** Our heuristic argument can be used to show the converse of the theorem, namely, if the column space of the design matrix $\mathsf{X}$ contains the eigenvector associated to the largest eigenvalue of the limiting covariance matrix, then the limiting variance of the BLUE $\hat{\beta}_\kappa$ is strictly positive. Indeed, note that the error distribution of $\eta_\kappa$ converges in distribution to the limiting distribution of $\eta_0$. In light of Fatou's lemma (see Ref.[10]), we have for all vectors $a$,

$$\liminf_{\kappa \longrightarrow 0} a^t\mathsf{V}(\hat{\beta}_\kappa)a \geq a^t\mathsf{V}(\hat{\beta}_0)a, \tag{21}$$

We may then use the heuristic to show that the right side of (21) strictly positive.

More general classes of limiting covariances matrix are discussed in Section 5.

# 4 Implications for sampling locations in experiments

The autocorrelation error model (15) of Section 3 provides a useful and simple framework in which to analyze parameter estimation from a large number of closely spaced measurements. In the context of magnetically confined plasmas, difficulty of access to the plasma typically implies that an increase in the number of measurements will be associated with a decrease in the spacing between measurements. As mentioned in the introduction, the noise in these devices arises in part from plasma turbulence, which may exhibit long-range characteristics. One may then be led to believe, in light of Section 3, that the increased correlations due to closer spacing could improve parameter estimation. We show below that this is not the case.

This section provides a detailed analysis of the variance of the BLUE for a single regression coefficient in the following setting: Suppose we observe the magnetic field at $n$ locations in the interval $[0, 1]$. For the purpose of our discussion, we take these points to be equidistant, that is $x_{n,i} = i/n$, $i = 0, \ldots, n$, (spacing $\Delta x = 1/n$). At each location $x_{n,i}$, we observe

$$Y_{n,i} = \mu + \eta_{n,i}, \quad i = 0, 1, \ldots, n, \tag{22}$$

where the disturbances have mean zero, variance $\mathbb{E}[\eta_{n,i}^2] = \sigma^2(x_{n,i}) = \tau^{-2}(x_{n,i})$ and correlation

$$\begin{aligned} \varrho_{ij} &= \tau(x_{n,i})\tau(x_{n,j})\mathbb{E}[\eta_{n,i}\eta_{n,j}] = \exp\left(-|x_{n,i} - x_{n,j}|/\delta\right) \\ &= \exp\left(-\frac{|i-j|}{\delta n}\right) \equiv \varrho^{|i-j|}. \end{aligned} \tag{23}$$

This is the autocorrelation model, equation (15) of Section 3, with $\varrho \equiv \exp(-(\delta n)^{-1})$. The parameter $\delta$ is interpreted as the correlation length, with $\delta \longrightarrow 0$ and $\delta \longrightarrow \infty$ corresponding to the uncorrelated error and fully correlated error models, respectively. Note that $\delta n = \delta/\Delta x$ is the ratio of the correlation length to the spacing between measurements. We shall further assume that $\tau$, a measure of the signal-to-noise ratio, is a smooth function of the sampling location.

In Section 3, we studied the limit of the variance $\mathcal{V}(n, \delta) \equiv V(\hat{\mu})$ as the correlation length $\delta$ approached infinity. In this section, we fix $\delta$ and study the behavior of the estimated mean $\hat{\mu}$ and its variance as the sample size $n$ goes to infinity. In this setting, as we increase the number of sampling locations within the unit interval, the correlation between neighboring measurements increases. While this is not the usual framework for asymptotic analysis, the result of this analysis can provide guidelines for the usefulness of acquiring additional data by increasing the number of measurements done for MHD equilibrium reconstructions, or other estimation problems where acquiring more data necessitates packing the measurements more closely in space or time.

THEOREM 2. *The inverse variance of the estimated mean for the autocorrelated model is, for $\delta n > 1$,*

$$\begin{aligned}\mathcal{V}^{-1}(n,\delta) &= \left[\frac{\delta}{2}\int_0^1 \tau'(s)^2 ds + \frac{1}{2\delta}\int_0^1 \tau(s)^2 ds\right]\left(1 + O((\delta n)^{-2})\right) \\ &\quad + \frac{1}{2}(\tau(0)^2 + \tau(1)^2) + \delta\, O(n^{-2}) + \delta^{-1}\, O(n^{-2}). \quad (24)\end{aligned}$$

**Remark.** This expression for the variance provides insight into the estimated mean. In the limit of very large correlations, the variance of the estimated mean is zero if $\int \tau'(s)^2 ds \neq 0$, i.e. if the signal-to-noise ratio varies over the measurement region. On the other hand, if $\int \tau'(s)^2 ds = 0$, the variance converges to $\sigma^2(0)$ for $\delta$ large. Indeed, vanishing of the integral of $\tau'(s)^2$ implies that $\tau(s) = \tau(0)$, in which case the variance converges with many measurements to

$$\lim_{n \to \infty} \mathcal{V}(n,\delta) = \frac{2\delta}{2\delta + 1}\sigma^2(0) \xrightarrow{\delta \to \infty} \sigma^2(0).$$

**Proof.** Let us denote $\tau_i \equiv \tau(x_{n,i})$, suppressing the $n$-dependence. The inverse of the covariance matrix of the autocorrelation model is

$$\Sigma^{-1} = \mathsf{S}^{-1}\mathsf{R}^{-1}\mathsf{S}^{-1},$$

with R as in equation (15),

$$\mathsf{R}^{-1} = \frac{1}{1-\varrho^2}\begin{bmatrix} 1 & -\varrho & 0 & \cdots & 0 \\ -\varrho & 1+\varrho^2 & -\varrho & \cdots & 0 \\ 0 & -\varrho & 1+\varrho^2 & \cdots & 0 \\ \vdots & \vdots & \vdots & \ddots & \vdots \\ 0 & 0 & 0 & \cdots & 1 \end{bmatrix},$$

and S the diagonal matrix with entries $\mathsf{S}_{ii} = \sigma_i = \tau_i^{-1}$. Using $\mathsf{X} = (1, 1, \ldots, 1)^t$ and equation (5), we find that the inverse of the variance $\mathcal{V}(n,\delta)$ of the estimated

mean is

$$\begin{aligned}
\mathcal{V}^{-1}(n,\delta) &= \sum_{i,j=0}^{n} (\mathsf{R}^{-1})_{ij} \tau_i \tau_j \\
&= \frac{1}{1-\varrho^2} \left[ \sum_{i=0}^{n} \tau_i^2 + \varrho^2 \sum_{i=1}^{n-1} \tau_i^2 - 2\varrho \sum_{i=0}^{n-1} \tau_i \tau_{i+1} \right] \\
&= \frac{1}{1-\varrho^2} \left[ (1-\varrho)^2 \sum_{i=0}^{n} \tau_i^2 + \varrho \sum_{i=0}^{n-1} (\tau_{i+1} - \tau_i)^2 \right. \\
&\quad \left. + \varrho(1-\varrho)(\tau(0)^2 + \tau(1)^2) \right].
\end{aligned} \qquad (25)$$

Since $x_{n,i} = i/n$, we can use the trapezoidal rule for numerical integration to approximate the sums (see Ref.[11])

$$\frac{1}{n} \sum_{i=0}^{n} \tau_i^2 = \int_0^1 \tau(s)^2 ds + \frac{1}{2n} \left( \tau(0)^2 + \tau(1)^2 \right) + O(n^{-2})$$

and, using the relation $(\tau_{i+1} - \tau_i) = (1/n)\tau_i' + (1/2n^2)\tau_i'' + O(n^{-3})$,

$$\begin{aligned}
n \sum_{i=0}^{n-1} (\tau_{i+1} - \tau_i)^2 &= \frac{1}{n} \sum_{i=0}^{n} (\tau_i')^2 - \frac{1}{n} \tau'(1)^2 + \frac{1}{n^2} \sum_{i=0}^{n} \tau_i' \tau_i'' + O(n^{-2}) \\
&= \int_0^1 \tau'(s)^2 ds - \frac{1}{2n} \left( \tau'(1)^2 - (\tau'(0))^2 \right) \\
&\quad + \frac{1}{2n} \int_0^1 \frac{d}{ds} (\tau'(s)^2) + O(n^{-2}) \\
&= \int_0^1 \tau'(s)^2 ds + O(n^{-2}).
\end{aligned}$$

Noting that $\varrho = \exp(-(\delta n)^{-1})$, it follows that

$$\begin{aligned}
\mathcal{V}^{-1}(n,\delta) &= \frac{1}{\delta} f\left((\delta n)^{-1}\right) \left( \int_0^1 \tau(s)^2 ds + O(n^{-2}) \right) + \frac{1}{2} (\tau(0)^2 + \tau(1)^2) \\
&\quad + \delta\, g\left((\delta n)^{-1}\right) \left( \int_0^1 \tau'(s)^2 ds + O(n^{-2}) \right),
\end{aligned}$$

where $f(x) = (1 - e^{-x})/x(1 + e^{-x})$ and $g(x) = xe^{-x}/(1 - e^{-2x})$. For $\delta n \gg 1$,

these functions both behave as

$$\lim_{\delta n \to \infty} f\left((\delta n)^{-1}\right) = \frac{1}{2} + O\left((\delta n)^{-2}\right)$$
$$\lim_{\delta n \to \infty} g\left((\delta n)^{-1}\right) = \frac{1}{2} + O\left((\delta n)^{-2}\right)$$

which produces the desired result, equation (24). □

Note that the usual result for uncorrelated errors ($\delta = 0$) may be recovered from (25) by simply setting $\varrho = 0$. In this case, the inverse variance is simply the sum of $n$ positive terms, implying the familiar relationship $\mathcal{V} \sim n^{-1}$ for $n$ uncorrelated measurements. In contrast, for $\delta n > 1$, there are no terms of order $n^{-1}$ in equation (24).

For a large number $n$ of highly correlated measurements ($\delta n \gg 1$), we may ignore the higher order terms in equation (24) of Theorem 2, and the variance converges to

$$\mathcal{V} = \frac{2\delta}{\int_0^1 \tau(s)^2 ds + \delta\left(\tau(0)^2 + \tau(1)^2\right) + \delta^2 \int_0^1 \tau'(s)^2 ds}. \tag{26}$$

Examples with $\delta > 0$ and linear inverse variance $\tau(s) = 1 + \alpha s$, with $\alpha = 1$ are shown in Fig. 3, with the results for finite $n$ summed numerically and the limit $n \to \infty$ from equation (26). The value of $\mathcal{V}$ converges rapidly as $n \to \infty$ except near $\delta = 0$, where $\mathcal{V} \sim n^{-1}$. The form of equation (26), including the behavior $\mathcal{V}(n, \delta) \sim 1/\delta$ for large $\delta$, is evident. The maximum of $\mathcal{V}(n, \delta)$ with respect to $\delta$ occurs at

$$\delta^2_{\mathcal{V}_{\max}} = \frac{\int_0^1 \tau(s)^2 ds}{\int_0^1 \tau'(s)^2 ds}.$$

That is, $\mathcal{V}$ decreases if the correlation length $\delta$ is larger than the typical scale for change of $\tau(x)$.

In Fig. 4 we show $\mathcal{V}(n, \delta)$ as a function of $n$ for three values of $\delta$. These results, similar to those of Ref.[6], show that $\mathcal{V}(n, \delta)$ converges to a positive value as $n \to \infty$, unless $\delta = 0$. There is an initial decrease, when $n \lesssim 1/\delta$; to the right of this region $\mathcal{V}(n, \delta)$ is nearly constant.

In Fig. 5 we show results for a constant inverse variance $\tau(s) = 1$ (i.e. $\alpha = 0$) both numerically for finite $n$ and the asymptotic result [equation (26)] for $n \to \infty$. Again, the results converge rapidly with $n$ except near $\delta = 0$, with $\mathcal{V} \to \tau^{-2} = 1$ as $\delta \to \infty$.

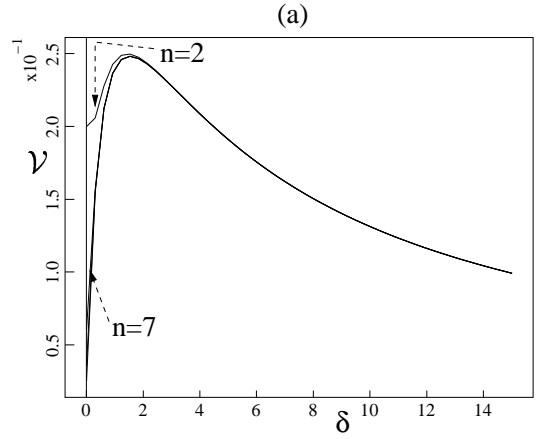

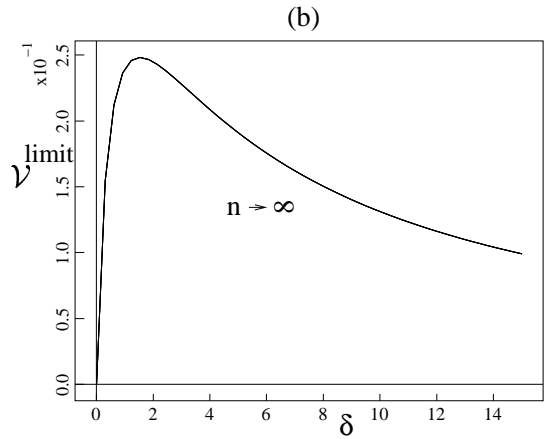

Figure 3: Estimate variance $\mathcal{V}$ as a function of $\delta$ for $\tau(1) = 1$ and $\alpha = 1$. The cases for $n = 2, 7$ are summed numerically and the case $\mathcal{V}^{\text{limit}}$ for $n \to \infty$ is from the analytic limit in Theorem 2.

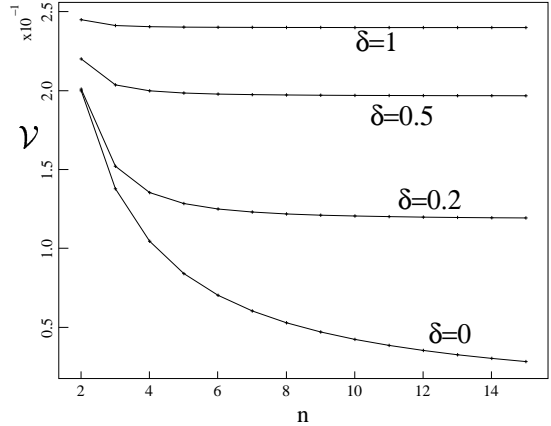

Figure 4: Estimate variance $\mathcal{V}$ as a function of $n$ for $\delta = 0, 0.2, 0.5, 1$, with $\alpha = 1$.

Notice that these results show that $\mathcal{V}(n, \delta)$ approaches a limiting curve as $n \to \infty$. Except for small $\delta$, the convergence is quite rapid due to the absence of corrections of order $n^{-1}$, as mentioned above. The results for $\alpha = 1$ show the generic situation of $\mathcal{V}(n, \delta) \to 0$ as $\delta \to \infty$; those with $\alpha = 0$ show the special situation in which $\mathcal{V}(n, \delta)$ approaches a positive constant in that limit. As long as $\delta > 0$, the estimate variance $\mathcal{V}(n, \delta)$ becomes constant with respect to $n$ for $n \gtrsim 1/\delta$ and has a finite limit as $n \to \infty$, showing that there is no advantage to be gained by increasing the number of measurements at points $x_i$ past $n \sim 1/\delta$.

## 5 Extensions and other considerations

In this section we discuss a few extensions of the above analysis.

### 5.1 Rank of $\tilde{\mathsf{X}}_{-1}$ less than $m$

In Section 3 we discussed the fact that for $\kappa \to 0$ the rank $r$ of the reduced design matrix $\tilde{\mathsf{X}}_{-1}$ is generically $m$, that in this case the vector $(e_1\sigma_1, ..., e_n\sigma_n)^t$ is not in the range of the original design matrix $\mathsf{X}$, and that the covariance of the estimate $\hat{\beta}_\kappa$ has $\text{trace}(\mathsf{V}(\hat{\beta}_\kappa)) \to 0$. On the other hand, if $r$, which is also the rank of $\tilde{\mathsf{X}}_{-1}^t \Lambda_{-1}^{-1} \tilde{\mathsf{X}}_{-1}$, is less than $m$, the covariance of $\hat{\beta}_\kappa$ does not go to zero. The noise-free subspace then has dimension $r$, meaning that $r$ linearly independent combinations of $\hat{\beta}_\kappa$ are determined exactly. The remaining $m-r$ combinations are subject to noise. Note that this case occurs when the vector of signed deviations $(e_1\sigma_1, \ldots, e_n\sigma_n)^t$ is in the range of the design matrix $\mathsf{X}$.

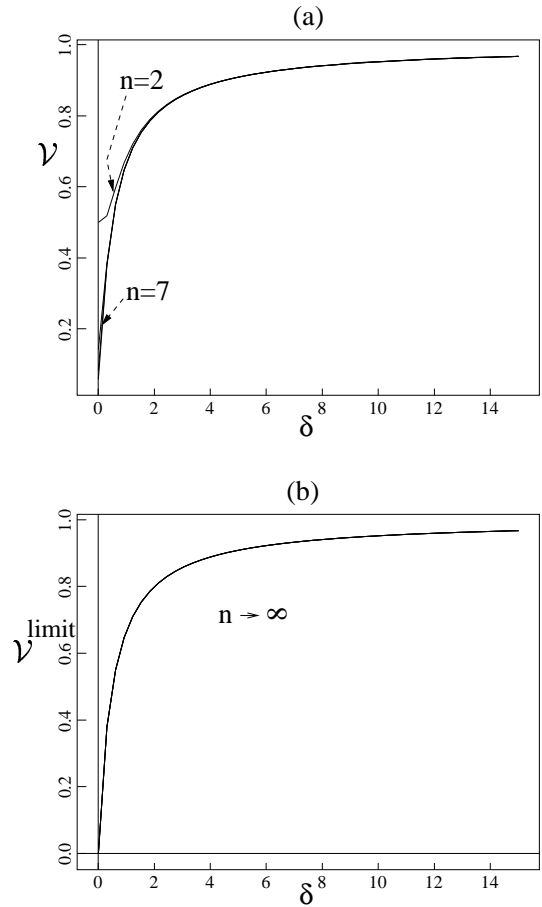

Figure 5: Estimate variance $\mathcal{V}$ as a function of $\delta$ for $\tau(1) = 1$, $\alpha = 0$ and (a) $n = 2, 7$, summed numerically, and (b) the limiting value $\mathcal{V}^{\text{limit}}$ for $n \to \infty$ from the analytic limit in Theorem 2.

## 5.2 Rank of $\Sigma$ approaches $r'$

In Sections 2 through 4, we have studied in depth the case in which the rank of $\Sigma$ approaches unity. To generalize, we suppose the rank goes in some limit to $r'$, with $1 \leq r' < n$. This can occur, for example, if two distinct and uncorrelated types of measurements are made. (This situation was present in the plasma reconstruction studies of Ref. [6], where magnetic field measurements external to the plasma and pressure measurements internal to the plasma were used.) For example, suppose one type of measurement, for $i = 1, ..., r'-1$ has a correlation matrix of the form $\mathsf{R}^{(1)} = \varrho_1^{|i-j|}$ and a second type, for $i = r', ..., n$ has $\mathsf{R}^{(2)} = \varrho_2^{|i-j|}$. We then have $\Sigma = \mathsf{SRS}$, with

$$\mathsf{R} = \left[\begin{array}{cc} \mathsf{R}^{(1)} & 0 \\ 0 & \mathsf{R}^{(2)} \end{array}\right].$$

Then for $\varrho_2 = 1$, but $|\varrho_1| < 1$, the rank of $\mathsf{R}$, and therefore the rank of $\Sigma$, equals $r'$.

In this case, the heuristic procedure described in Section 3 leads, in the limit $\varrho_2 \to 1$, to a linear system of equations subjected to noise

$$Z_j = v_j^t \mathsf{X}\beta + \lambda_j^{1/2}\xi_j \quad j = 1, ..., r',$$

and a noise-free subspace

$$Z_j = v_j^t \mathsf{X}\beta \quad j = r'+1, \cdots, n.$$

If $n - r' \geq m$, the second set of equations determines $\hat{\beta}$ exactly and the estimate variance is zero. If, on the other hand, $n - r' < m$, the second set of equations has a null space of dimension $m - n + r'$. That is, there are $n - r'$ linearly independent linear combinations of the $\beta_i$ that are determined exactly. In other words, the estimate covariance matrix has rank $n - r'$, i.e. $n - r'$ zero eigenvalues and $m - n + r'$ nonzero eigenvalues.

## 6 Summary and discussion

In its fundamental form, the main result of this paper, given in Section 3, is the following: in the limit of strong correlations characterized by a single correla-tion length $\delta \to \infty$, the covariance matrix V of the estimate $\hat{\beta} = (\hat{\beta}_1, \ldots, \hat{\beta}_m)$ generically vanishes. That is, its trace, the total variance of V, vanishes. The exceptions to this rule occur when the vector of signed deviations $(e_1\sigma_1, ..., e_n\sigma_n)^t$ of the measurements is in the range space of the design matrix X. We explained the decrease and eventual vanishing of trace(V) by means of a simple example in Section 2, and also showed the relationship between this phenomenon and negative weighting, in which the estimate is a weighted average of the measurements, with some weights negative.

This result is so surprising that it suggests a "free lunch" possibility. The idea that stronger correlations can be obtained simply by packing in closer measurements has been studied in Section 5. It is found that the covariance of the estimate does indeed decrease as the number $n$ of measurements increases, but this decrease flattens when $\delta \gtrsim \Delta x$, where $\Delta x$ is the spacing between measurements. The interpretation of this result is that for the variance to decrease with increasing number of measurements, the measurement spacing must not be much smaller than the correlation length. Further, from (26) we have concluded that the variance decreases with correlation length $\delta$ if $\delta$ is greater than the typical length scale for variations in the signal-to-noise ratio.

We have addressed other related issues in Section 6. The first is the situation in which the reduced design matrix is not of full rank. In this case the noise-free subspace has dimension $r < m$, meaning that $m - r$ linearly independent combinations of $\hat{\beta}_i$ are determined exactly, and the other $r$ combinations are not determined exactly and therefore involve the noise. We also considered a generalization of the condition related to characterization of the correlations by a single parameter $\kappa$, i.e. the condition that, as $\kappa \to 0$, the rank of the data covariance matrix $\Sigma$ goes from $n$ to unity. The generalization deals with cases in which this rank decreases from $n$ to $r'$ as some parameter is varied. In this case, the result is unchanged if $r' \geq m$; if $r' < m$, however, the estimate is not completely determined in the limit of large correlations, but a linearly independent set of $n - r'$ combinations are determined exactly.